\documentclass{amsart}
\usepackage{amsmath, amsthm, latexsym, amssymb, amsfonts, color, epsfig}

\newenvironment{pf}{\proof[\proofname]}{\endproof}
\theoremstyle{plain}
\newtheorem{Th}{Theorem}[section]
\newtheorem{Cor}[Th]{Corollary}
\newtheorem{Prop}[Th]{Proposition}
\newtheorem{Lemma}[Th]{Lemma}
\numberwithin{equation}{section} \theoremstyle{definition}
\newtheorem{Rem}[Th]{Remark}

\newtheorem{Example}{Example}
\newtheorem{Def}[Th]{Definition}

\newcommand{\cal}[1]{\mathcal{#1}}
\newcommand{\C}{\mathbb C}
\newcommand{\Q}{\mathbb Q}
\newcommand{\Z}{\mathbb Z}
\newcommand{\R}{\mathbb R}

\newcommand{\D}{\Delta}

\newcommand{\cF}{\cal F}

\newcommand{\cO}{\cal O}
\newcommand{\cU}{\cal U}
\newcommand{\cX}{\cal X}
\newcommand{\cY}{\cal Y}
\newcommand{\cZ}{\cal Z}

\newcommand{\p}{\partial}
\newcommand{\T}{{\mathbb T}}

\newcommand{\sig}{\sigma}
\newcommand{\Sig}{\Sigma}

\newcommand{\Res}{\operatorname{Res}}

\newcommand{\sgn}{\operatorname{sgn}}

\newcommand{\cdeg}{\operatorname{c.\!deg}}

\newcommand{\Tr}{\operatorname{Tr}}

\newcommand{\Hom}{\operatorname{Hom}}

\newcommand{\rs}[1]{Section~\ref{S:#1}}

\newcommand{\rp}[1]{Proposition~\ref{P:#1}}

\newcommand{\rex}[1]{Example~\ref{ex:#1}}
\newcommand{\re}[1]{(\ref{e:#1})}
\newcommand{\rc}[1]{Corollary~\ref{C:#1}}
\newcommand{\rt}[1] {Theorem~\ref{T:#1}}

\newcommand{\rf}[1]{Figure~\ref{F:#1}}

\begin{document}


\title{Toric residue and combinatorial degree}
\author{Ivan Soprounov}
\address{Department of Mathematics and Statistics,
University of Massachusetts, Amherst, MA 01003}
\email{isoprou@math.umass.edu}

\begin{abstract} Consider an $n$-dimensional projective toric variety $X$
defined by a convex lattice polytope $P$. David Cox introduced the
toric residue map given by a collection of $n+1$ divisors
$(Z_0,\dots,Z_n)$ on $X$. In the case when the $Z_i$ are $\T$-invariant
divisors whose sum is $X\setminus\T$ the toric residue map
is the multiplication by an integer number. We show that this
number is the degree of a certain map from the boundary of the
polytope $P$ to the boundary of a simplex. This degree can be computed
combinatorially. We also study radical monomial ideals $I$ of the homogeneous
coordinate ring of $X$. We give a necessary and sufficient condition
for a homogeneous polynomial of semiample degree to belong to $I$
in terms of geometry of toric varieties and combinatorics of fans.
Both results have applications to the problem of constructing an element
of residue one for semiample degrees. 

\end{abstract}

\maketitle

\section*{Introduction}

Toric residue is defined for every collection of $n+1$
divisors on an $n$-dimensional complete toric variety as long as
they do not have a common point. It appears in a variety of contexts, e.g. 
in mirror symmetry \cite{BM} or in sparse polynomial systems \cite{CaD}, \cite{CDS}.
Toric residues were first introduced by D.~Cox
in \cite{Coxres} in the case when the divisors are ample and
linearly equivalent. In \cite{CCD} Cattani, Cox, and Dickenstein
extended the definition to the general case and revealed a
connection between the toric residue and the sum of local
Grothendieck residues in the torus. They also provided an
algorithm for computing the toric residue when the divisors are
ample using Gr\"obner bases. Another approach
was taken by D'Andrea and Khetan in \cite{AK} where they compute
the toric residue as a quotient of two determinants. One of the
key ideas in both methods is reduction to a particular choice of
$n+1$ $\T$-invariant divisors for which the toric residue is one.
This motivates the problem of computing the toric residue for an
arbitrary choice of $n+1$ $\T$-invariant divisors.

Let $X$ be a projective toric variety of dimension $n$ defined by a convex lattice
polytope $P$. Consider $n+1$ effective $\T$-invariant divisors
$Z_0,\dots,Z_n$ on $X$ whose sum $Z$ equals $X\setminus\T$ (or more
generally the support of $Z$ equals $X\setminus\T$).
Then the toric residue map for
$Z_0,\dots,Z_n$ has a nice combinatorial description: Notice that
the irreducible components of the $Z_i$ correspond to the facets
of the polytope $P$. Thus, we get a ``coloring'' of the facets of
$P$ into $n+1$ colors: a facet has color~$i$ if the corresponding
irreducible divisor appears in $Z_i$. Let us now pick a color
$0\leq i\leq n$ for each facet of the $n+1$ facets of the standard
$n$-simplex $\D$, so that all facets have different colors.
Choose a continuous piecewise linear map $f:\p P\to \p\D$ that
matches the colors. Such a map exists as long as $Z_0\cap\dots\cap
Z_n$ is empty. The degree of $f$ is independent
of the choice of $f$ and is called the combinatorial degree of the
coloring of $P$ (we give the precise definition in \rs{cdeg}). We
prove that the toric residue for $Z_0,\dots,Z_n$ coincides with
the combinatorial degree of the coloring defined by the $Z_i$ (\rt{main1},
\rt{main2}). The combinatorial degree can be computed explicitly
as a signed number of certain complete flags of faces of~$P$
(see \rt{flags}).

In the second part of the paper we consider a radical monomial ideal
$I$ of the homogeneous coordinate ring
$S=\bigoplus_{\alpha\in A_{n-1}(X)}S_\alpha$ of $X$ (see \cite{Coxhom}).
Let $z_1,\dots,z_m$ be the minimal generating set of $I$.
Consider a homogeneous polynomial $F\in S_\alpha$ of semiample degree
$\alpha$ (i.e. the
line bundles corresponding to $\alpha$ are generated by global sections).
It follows from \cite{CCD}, Section~2 that if $\alpha$ is $\Q$-ample
the polynomial~$F$ belongs to $I$ as long as $Z_1\cap\dots\cap Z_m$
is empty, where $Z_i$ is the zero set of $z_i$ on~$X$. In \rt{monomial}
we extend this observation.
It is known that if $\alpha$ is semiample there exists
a unique toric variety $X'$ and a surjective morphism $\pi:X\to X'$ such that
$\alpha$ is the pull-back of an ample class on $X'$ (see \cite{M},
Theorem~1.2). We show that a generic
polynomial $F$ of semiample degree lies in $I$ if and only if
$\pi(Z_1)\cap\dots\cap\pi(Z_m)$ is empty on $X'$. This condition also
has a combinatorial interpretation (\rt{monomial}, part (4)).

Let us remark that it is still an open question how to compute the toric
residue map for arbitrary non-ample divisors even when we know
that the map is an isomorphism (e.g. for big and nef divisors).
In the final section we apply our results to this problem. We obtain
a combinatorial condition when it is possible to construct an element
of toric residue one. More results in this direction are obtained in the joint work with
A.~Khetan \cite{KS}.

{\bf Acknowledgments.} I would like to thank Eduardo Cattani for
introducing me to the subject, his constant attention to this work,
and numerous discussions. I am also grateful to David Cox for
helpful discussions, useful comments, and a
suggestion on how to improve the result of
\rt{monomial}. The result of \rt{main1} was presented at
the MEGA 2003 conference in Kaiserslautern.

\section{Combinatorial Degree}\label{S:cdeg}

\subsection{Polyhedral sets}
A polyhedral set $X$ is a finite union of convex compact polyhedra
intersecting in faces. We assume that they all are embedded in
some Euclidean space $E$, thus $X$ is a topological space with
topology inherited from $E$. The dimension of $X$ is the maximum
of dimensions of the polyhedra it contains.

\begin{Example}\label{ex:boundary} {\it Boundary subdivision.}
Let $P\subset E$ be a convex $n$-dimensional polytope.
Then the boundary $\partial P$ is a polyhedral set of
dimension~$n-1$. More generally, any polyhedral subdivision of $\partial P$
is a polyhedral set of dimension $n-1$.
\end{Example}

For the purposes of this paper it is enough
to consider only polyhedral sets of \rex{boundary}.

Let $X$ be a polyhedral set and $\cF(X)$ the set of all faces of
all polyhedra appearing in $X$. The set $\cF(X)$ is a finite
partially ordered set by inclusion.

\begin{Def}
Let $X$, $Y$ be two polyhedral sets and $\psi:\cF(X)\to\cF(Y)$ a
map that preserves the partial ordering. A continuous map
$f_{\psi}:X\to Y$ is called a {\it characteristic map for $\psi$} if the
image of every face $G\in\cF(X)$ under $f_{\psi}$ lies in
$\psi(G)$.
\end{Def}

For any map $\psi:\cF(X)\to\cF(Y)$ there exists a characteristic map
$f_{\psi}:X\to Y$. One can construct a
piecewise linear $f_\psi$ using barycentric subdivisions of~$X$
and~$Y$ (see \cite{S}, Proposition 2.1).

It is easy to see that if $f_{\psi}$ and $f'_{\psi}$ are two
characteristic maps for $\psi$ then
for any $0\leq t\leq 1$ the map $f^t_{\psi}=(1-t)f_{\psi}+tf'_{\psi}$
is also characteristic for $\psi$.

\begin{Def} Let $X$ and $Y$ be $n$-dimensional polyhedral sets
and $\psi:\cF(X)\to\cF(Y)$ a map that preserves the partial
ordering. Let $f_{\psi}$ be any characteristic map for $\psi$. The
induced map of the $n$-th homology groups
\begin{equation}\label{e:map}
H_n(f_{\psi}):H_n(X)\to H_n(Y)
\end{equation}
is called the {\it combinatorial degree map} of $\psi$.
\end{Def}

By the above the combinatorial degree map
is independent of the choice of a characteristic map $f_{\psi}$.
We denote it by $\cdeg\psi$.

Next we will see that the combinatorial degree map is invariant under
subdivisions of~$X$.

\begin{Def} A polyhedral set $X'$ is called a {\it subdivision} of a polyhedral
set $X$ if every polyhedron in $X$ is a finite union of polyhedra
in $X'$, and the support of $X$ equals the support of $X'$.
We say that $\psi':\cF(X')\to\cF(Y)$ {\it refines}
$\psi:\cF(X)\to\cF(Y)$ if $G'\subset G$ implies
$\psi'(G')\subset\psi(G)$ for every $G'\in\cF(X')$ and
$G\in\cF(X)$.
\end{Def}

\begin{Prop}\label{P:subdivision} Let $X$ and $Y$ be $n$-dimensional
polyhedral sets. Let $X'$ be a subdivision of $X$. Consider
$\psi:\cF(X)\to\cF(Y)$ and
$\psi':\cF(X')\to\cF(Y)$ that
refines~$\psi$. Then $\cdeg\psi'=\cdeg\psi$.
\end{Prop}

\begin{pf} Let $f_{\psi'}:X'\to Y$ be any characteristic map for $\psi'$.
Since $\psi'$ refines $\psi$ the map $f_{\psi'}$ considered as a
map from $X$ to $Y$ is characteristic for $\psi$. Therefore
$\cdeg\psi'=H_n(f_{\psi'})=\cdeg\psi$.
\end{pf}

\subsection{Combinatorial degree for polytopes}

Now assume that $X$ and $Y$ are boundary subdivisions of some convex
$n$-dimensional polytopes $P$ and $Q$ as in \rex{boundary}.
In this case we have isomorphisms $H_{n-1}(X)\cong\Z$ and $H_{n-1}(Y)\cong\Z$
given by a choice of orientations of $P$ and $Q$. The map \re{map} is then the
multiplication by an integer number. This number will be
called the {\it combinatorial degree} of $\psi$ and denoted
also by $\cdeg\psi$.

\begin{Rem} The combinatorial degree is a global analog of combinatorial
coefficients, first introduced by O.~Gelfond and A.~Khovanskii
in \cite{GKh}. Combinatorial coefficients appear in the Gelfond-Khovanskii
residue formula and the Khovanskii product of roots formula for
a polynomial system whose Newton polytopes have generic relative
position (see \cite{GKh2} and \cite{Kh} for details).
\end{Rem}

The next theorem shows how to compute the combinatorial degree
for polytopes. We will need the following definition.

\begin{Def} Let $X$ be a boundary subdivision of an
$n$-dimensional convex oriented polytope $P$. Consider a complete
flag in $X$ (a maximal chain of elements of $\cF(X)$):
$$\cX:\ X_0\subset\dots\subset X_{n-1},\quad \dim X_i=i.$$
For every $1\leq i\leq n$ choose a vector $e_i$ which begins at $X_0$
and points strictly inside~$X_i$, where $X_n=P$.
Define the {\it sign} of the flag $\cX$ to be 1
if $(e_1,\dots,e_n)$ gives a positive oriented frame for $P$,
and $-1$ otherwise.
\end{Def}

Let $X$ and $Y$ be boundary subdivisions of $n$-dimensional convex
oriented polytopes $P$ and $Q$, and $\psi:\cF(X)\to\cF(Y)$. For
complete flags $\cX:\, X_0\subset\dots\subset X_{n-1}$ in~$X$ and
$\cY:\ Y_0\subset\dots\subset Y_{n-1}$ in $Y$ we will write
$\psi(\cX)=\cY$ if and only if $\psi(X_i)=Y_i$ for all $0\leq
i\leq n-1$.

\begin{Th}\cite{S}\label{T:flags}
Fix any complete flag $\cY$ in $Y$.
Then the combinatorial degree of $\psi$ is equal
to the sign of $\cY$ times
the number of all complete flags $\cX$ in $X$
counted with signs, such that $\psi(\cX)=\cY$:
$$\cdeg\psi=\sgn\cY\sum_{\psi(\cX)=\cY}\sgn\cX.$$
\end{Th}


Now let $P$ be an $n$-dimensional convex polytope in $E$, $\dim E=n$, and
consider its {\it polar} polytope $P^\circ$ in the dual space $E^*$.
Let $X=\p P$ and $X^*=\p(P^\circ)$.
Recall that there is a one-to-one and order-reversing
correspondence $G\mapsto G^*$ between $\cF(X)$ and $\cF(X^*)$
such that $\dim G^*=n-1-\dim G$.

\begin{Def} Let $X=\p P$, $Y=\p Q$ for convex $n$-dimensional polytopes $P$ and~$Q$
in $E$.  Let $\psi:\cF(X)\to\cF(Y)$ be a map of partially
ordered sets. Define the {\it dual map}
$\psi^*:\cF(X^*)\to\cF(Y^*)$ by putting $\psi^*(G^*)=(\psi(G))^*$
for any $G\in\cF(X)$.
\end{Def}

It is readily seen that $\psi^*$ preserves the partial order.

\begin{Prop}\label{P:dual} Let $P$, $Q$ be convex $n$-dimensional polytopes in
$E$, and $X=\p P$, $Y=\p Q$. Consider any $\psi:\cF(X)\to\cF(Y)$
and its dual $\psi^*:\cF(X^*)\to\cF(Y^*)$. Then
$\cdeg\psi=\cdeg\psi^*$.
\end{Prop}

\begin{pf} This follows from \rt{flags}. Indeed, every complete
flag $\cX$ in $X$ corresponds to a unique complete flag $\cX^*$ in
$X^*$ and $\sgn\cX=\lambda\sgn\cX^*$, where $\lambda=\pm 1$ and is
the same for all $\cX$.\footnote{In fact, $\lambda=(-1)^{\frac{n(n+1)}{2}}$.}
Then for any fixed flag $\cY$ in $Y$
$$\cdeg\psi=\sgn\cY\sum_{\psi(\cX)=\cY}\sgn\cX=
\lambda\sgn\cY^*\sum_{\psi^*(\cX^*)=\cY^*}\lambda\sgn\cX^*=\cdeg\psi^*.$$
\end{pf}

\subsection{Simplicial coloring}

Consider a convex $n$-dimensional polytope $P$.

\begin{Def}
Let $C_0,\dots,C_n$ be closed subsets of $\p P$, each set being the
union of some facets of $P$. We say that $C=(C_0,\dots,C_n)$ forms
a {\it simplicial coloring} of (the boundary of) $P$ if
\begin{enumerate}
\item $C_0\cup\dots\cup C_n=\partial P$,
\item $C_0\cap\dots\cap C_n=\emptyset$.
\end{enumerate}
\end{Def}

Let $\D$ denote the standard $n$-simplex:
$$\D=\{y=(y_0,\dots,y_n)\in\R^{n+1}\ |\ y_0+\dots+y_n=1,\, 0\leq y_i\leq 1\}.$$
For each $1\leq k\leq n$ the codimension $k$ faces of $\D$
are $$\D_{i_1\dots i_k}=\{y\in\D\ |\ y_{i_1}=\dots=y_{i_k}=0\},\quad
\text{where }\ 0\leq i_1<\dots<i_k\leq n.$$
Then every simplicial coloring $C=(C_0,\dots,C_n)$
defines a map
$$\psi_C: \cF(\partial P)\to\cF(\partial\D),$$
given by $\psi_C(G)=\D_{i_1\dots i_k}$ if $G$ belongs to
$C_{i_1}\cap\dots\cap C_{i_k}$ with $k$ maximal. Clearly $\psi_C$
respects the partial order.

\begin{Def} Fix orientations of the polytope $P$ and the simplex $\D$.
The {\it combinatorial degree of a simplicial coloring}
$C=(C_0,\dots,C_n)$ of $P$ is the combinatorial degree $\cdeg\psi_C$ of the
corresponding map $\psi_C: \cF(\partial P)\to\cF(\partial\D).$
\end{Def}

It follows from the definition that $\cdeg\psi_C$ is alternating on
the order of the $C_i$.

Recall that the {\it normal fan} $\Sig_P$ of a polytope $P$ in $E$
is a complete fan in $E^*$ whose cones are
\begin{equation}\label{e:normalfan}
\sig_{G}=
\{v\in E^*\ |\ \langle u,v \rangle\geq\langle u',v \rangle, \text{
for all } u\in P,\ u'\in G \},
\end{equation}
for every face $G$ of $P$. Note
that the 1-dimensional cones (rays) of $\Sig_P$ are generated by
the inner normals to the facets of $P$. Thus, every coloring of
$P$ defines a ``coloring'' of the rays of $\Sig_P$ into $n+1$
colors. We arrive at the following definition.

\begin{Def} Let $\Sig$ be a fan in $\R^n$. A decomposition
of the set $\Sig(1)$ of rays of $\Sig$ into a union of
$n+1$ subsets
$$\Sig(1)=\Lambda_0\cup\dots\cup\Lambda_n$$
is called a {\it coloring of $\Sig$}. We say a coloring is {\it
disjoint} if the union above is disjoint.

A coloring is called {\it simplicial} if no (maximal) cone of
$\Sig$ contains rays of all $n+1$ colors, i.e. every set of rays
$\{\rho_{i}\in\Lambda_i,\,0\leq i\leq n\}$ is not contained in any
of the (maximal) cones of $\Sig$.
\end{Def}

Note that we allow multiple colors of a ray unless the coloring is
disjoint.

Clearly, every simplicial coloring of a polytope $P$ defines a
simplicial coloring of its normal fan $\Sig_P$. Conversely, if
$\Sig$ is a {\it projective} fan (i.e. the normal fan of a polytope)
then every
simplicial coloring of $\Sig$ defines a simplicial coloring
of any polytope whose normal fan is $\Sig$. It is not hard
to see that these colorings will all have the same
combinatorial degree (for example one can deduce it from \rt{flags}).
We will call it the {\it combinatorial
degree of the simplicial coloring of~$\Sig$}.

Applying \rt{flags} we obtain an explicit formula
for the combinatorial degree of a coloring of a fan $\Sig$ in terms
of complete flags of faces of $\Sig$. For simplicity and
for further purposes we will state this formula in the case
when $\Sig$ is {\it simplicial} (i.e. all cones of $\Sig$ are simplicial)
and the coloring is disjoint.

Fix $0\leq i_1<\dots<i_k\leq n$. Then we say that a cone $\sig$ of
$\Sig$ is {\it $(i_1,\dots,i_k)$-colored} if for every $1\leq
j\leq k$ the cone $\sig$ contains a ray from $\Lambda_{i_j}$, and $k$ is maximal.

Let $\sig$ be a maximal simplicial $(i_1,\dots,i_n)$-colored cone
for a simplicial disjoint coloring of $\Sig$.
We say it has {\it positive orientation} if a
collection of generators $(e_1,\dots,e_n)$, where $e_j$ generates
the ray from $\Lambda_{i_j}$, gives a positive oriented frame for
$\R^n$.

\begin{Cor}\label{C:simplicial}
Let $\Sig$ be a projective simplicial fan in $\R^n$ and
$$\Sig(1)=\Lambda_0\sqcup\dots\sqcup\Lambda_n$$
a disjoint simplicial coloring. Fix any $k$, $0\leq k\leq n$. Then the
combinatorial degree of the coloring is equal to $(-1)^k$ times
the number of maximal $(0,\dots,\widehat k,\dots, n)$-colored
cones, counted with orientations.
\end{Cor}

Now suppose that $\Sig'$ is a refinement of $\Sig$ with the same
set of rays. Then any simplicial coloring of $\Sig$ is also a
simplicial coloring of $\Sig'$. The following theorem shows that
these colorings will have the same combinatorial degree.

\begin{Th}\label{T:maindegree} Let $\Sig$ and $\Sig'$ be projective
fans such that
$\Sig'\to\Sig$ is a refinement and $\Sig(1)=\Sig'(1)$. Consider a
simplicial coloring of $\Sig$. The induced coloring of $\Sig'$ is
simplicial and has the same combinatorial degree.
\end{Th}

\begin{pf} The fact that the induced coloring of $\Sig'$ is
simplicial follows by definition.

Let $P$ and $P'$ be polytopes whose normal fans are $\Sig$ and
$\Sig'$, respectively. Furthermore, let $C$ (resp. $C'$) be the
simplicial coloring of $P$ (resp. $P'$) defined by the simplicial
coloring of $\Sig$ (resp. $\Sig'$). We need to show that
$\cdeg\psi_C=\cdeg\psi_{C'}$.

We can assume that the origin lies in the interior of $P$ and $P'$
and so the fans $\Sig$ and $\Sig'$ consist of the cones over the
proper faces of $P^\circ$ and $(P')^\circ$, respectively. By
\rp{dual} it is enough to show that
$\cdeg\psi_C^*=\cdeg\psi_{C'}^*$ for the dual maps.

Let $X$ denote the boundary of $P^\circ$.
Projecting the proper faces of $(P')^\circ$ onto the
faces of $X$ along the rays of $\Sig$ we get a subdivision $X'$
of $X$. It is clear that we get an order-preserving
bijection between $\cF(X')$ and the proper faces of $(P')^\circ$,
and that the composition map $\psi':\cF(X')\to\cF(\p(\D^\circ))$
has the same combinatorial degree as $\psi_{C'}^*$.

It remains to show that $\psi':\cF(X')\to\cF(\p(\D^\circ))$ refines
$\psi_C^*: \cF(X)\to\cF(\p(\D^\circ))$.
Indeed, $\psi_C^*$ sends an element $G\in\cF(X)$ to $\D_{i_1\dots i_k}^*$
if and only if the cone over $G$ is $(i_1,\dots,i_k)$-colored.
Suppose $G'\in\cF(X')$ is contained in $G$. Then the rays of the
cone over $G'$ are also rays of the cone over $G$ and, hence,
will have the same colors $(i_1,\dots,i_k)$ (all of them or
fewer). Therefore, $\psi'(G')\subset\D_{i_1\dots
i_k}^*$. By \rp{subdivision} $\cdeg\psi'=\cdeg\psi_{C}^*$
which completes the proof.
\end{pf}

\section{Toric residue}

Let $X=X(\Sig)$ be an $n$-dimensional complete toric variety
determined by a complete rational fan $\Sig$ (see for example
\cite{F}). Following the notation of \cite{F} we let $N$ denote a
lattice of rank $n$ and $N_{\R}$ the real vector space
$N\otimes\R$ which contains the fan $\Sig$. Also let
$M=\Hom(N,\Z)$ be the dual lattice and $M_{\R}$ the corresponding
dual space $M\otimes\R$.

Let $\Sig(1)=\{\rho_1,\dots,\rho_r\}$ be the set of rays of
$\Sig$. Each ray $\rho_i\in\Sig(1)$ determines a $\T$-invariant
irreducible divisor $D_i$ on $X$. The variety $X$ has the
homogeneous coordinate ring $S=\C[x_1,\dots,x_r]$ graded by the
Chow group $A_{n-1}(X)$ so that a monomial
$x^a=\prod_{j=1}^{r}x_j^{a_j}$ has degree
$\deg(x^a)=[\sum_{j=1}^{r}a_{j}D_{j}]\in A_{n-1}(X)$ (see \cite{Coxhom}).
Denote by $S_\alpha$ the graded piece of $S$ consisting of all
polynomials of degree $\alpha\in A_{n-1}(X)$. As shown in \cite{Coxhom} $S_\alpha$
is canonically isomorphic to the global sections of the sheaf $\cO_X(D)$ on $X$,
where $\alpha=[D]$.

We will recall the definition of the toric residue \cite{Coxres, CCD}.
Consider $n+1$ homogeneous polynomials $F_i\in S_{\alpha_i}$, for $0\leq i\leq n$.
Their critical degree is defined to be
$$\rho=\sum_{i=0}^n\alpha_i-\sum_{j=1}^r\deg(x_j).$$
Then every polynomial $H$ of degree $\rho$ defines a meromorphic
$n$-form on $X$:
$$\omega_F(H)=\frac{H\Omega}{F_0\dots F_n},$$
where $\Omega$ is the Euler form \cite{CCD}. We use $F$ to denote
the vector $(F_0,\dots,F_n)$.

Suppose that $F_i$ do not vanish simultaneously on $X$. Then $X$
has an open cover $\cU=\{U_i,\,0\leq i\leq n\}$, where $U_i=\{x\in X\ |\
F_i(x)\neq 0\}$. Therefore $\omega_F(H)$ defines a \v Cech
cohomology class $[\omega_F(H)]\in H^n(X,\widehat\Omega_X^n)$,
relative to the cover $\cU$. Here $\widehat\Omega_X^n$ denotes the
sheaf of Zariski $n$-forms on $X$. This class $[\omega_F(H)]$ is
alternating on the order of the $F_i$ and is zero if $H$ belongs
to the ideal $\langle F_0,\dots,F_n\rangle$. The {\it toric
residue map}
$$\Res_F:S_\rho/\langle F_0,\dots,F_n\rangle\to\C,$$
is given by $\Res_F(H)=\Tr_X([\omega_F(H)])$, where $\Tr_X$ is the
trace map on $X$.

As proved in \cite{CCD}, when $X$ is simplicial the toric
residue is in fact the sum of local Grothendieck residues:

\begin{Th}\label{T:local}\cite{CCD} Let $X$ be a complete simplicial toric
variety. Let $F_0,\dots,F_n$ be homogeneous polynomials
which do not vanish simultaneously on $X$ and suppose the set
$Z_{\hat{k}}=\{x\in X\ |\ F_i(x)=0,\ 0\leq i\leq n\,,i\neq k\}$ is finite,
for some $k$. Then for any homogeneous polynomial $H$ of the critical degree $\rho$
$$\Res_F(H)=(-1)^k\sum_{x\in Z_{\hat k}}
\Res_{x}\left(\frac{(H/F_k)\Omega}{F_0\dots\widehat F_k\dots F_n}\right).$$
\end{Th}

We consider a special case when $F_i=z_i$ are monomials whose
product is the product of the variables, $z_0\dots z_n=x_1\dots
x_r$. In this case the critical degree $\rho$ is zero and the
toric residue is the multiplication by $\Res_{z}(1)$.

On the other hand, since the variables $x_j$ correspond to the
rays $\rho_j$ of $\Sig$, the monomials $z_0,\dots,z_n$ define a
disjoint coloring of $\Sig$:
$$\Sig(1)=\Lambda_0\sqcup\dots\sqcup\Lambda_n,$$
where the ray $\rho_j$ lies in $\Lambda_i$ if and only if $x_j$
divides $z_i$. Clearly, this coloring is simplicial if and only if
$z_0,\dots,z_n$ do not vanish simultaneously on $X$.

In the following theorem we show that $\Res_{z}(1)$ is equal to
the combinatorial degree of the coloring defined by
$z_0,\dots,z_n$.

\begin{Th}\label{T:main1} Let $X$ be an $n$-dimensional projective toric
variety defined by a projective fan $\Sig$. Let $z_0,\dots,z_n$ be
monomials in the homogeneous coordinate ring $S=\C[x_1,\dots,x_r]$
such that
\begin{enumerate}
\item $z_0\dots z_n=x_1\dots x_r$,
\item $z_0,\dots,z_n$ do not vanish simultaneously on $X$.
\end{enumerate}
Then the toric residue $\Res_{z}(1)$ equals the combinatorial
degree of the simplicial coloring of $\Sig$ defined by
$z_0,\dots,z_n$.
\end{Th}

\begin{pf}
First assume that $X$ is simplicial, and so $\Sig$ is a simplicial
fan. Let $Z_i$ be the zero locus of $z_i$ on $X$.
Suppose the intersection $Z_{\hat 0}=Z_1\cap\dots\cap Z_n$ is
infinite. Since it is a union of orbits it must contain a 1-dimensional orbit
$O_\tau$, for some cone~$\tau$ of dimension $n-1$.
But this would imply that there are
$n$ distinct irreducible components $Z_i'\subset Z_i$, $1\leq i\leq n$
which contain $O_\tau$, i.e. the cone
$\tau$ contains $n$ distinct rays of~$\Sig$, which is impossible
when $\Sig$ is simplicial. Therefore, $Z_{\hat 0}$ is finite.

By \rt{local} the toric residue is the sum of the Grothendieck residues over
the points of $Z_{\hat 0}$. But these points are the closed orbits
$O_\sig$ of $X$ that correspond to $(1,\dots,n)$-colored maximal
cones $\sig$ of $\Sig$. Therefore by \rc{simplicial} one only
needs to check that the Grothendieck residue at $O_\sig$ is equal
to $\pm 1$ depending on the orientation of $\sig$, which is
straightforward.

In the  general case let $X'=X(\Sig')$ be the toric variety
determined by a simplicial refinement $\Sig'$ of $\Sig$ with the
same set of rays. Then $X$ and $X'$ have the same homogeneous
coordinate ring and the birational morphism $X'\to X$ maps every
$Z'_i$ defined by $z_i$ on $X'$ onto $Z_i$. By the functorial
property of the trace map the toric residues $\Res_{z}(1)$ on $X$
and on $X'$ are equal. On the other hand, by
\rt{maindegree} the combinatorial degrees of the colorings of
$\Sig$ and $\Sig'$ are also equal, and the theorem follows.
\end{pf}

We will now consider a more general case when $z_i$ are
any (monic) monomials whose product is divisible by the
product of the variables, $x_1\dots x_r\,|\,z_0\dots z_n$. In this
case the quotient $z_0\dots z_n/{x_1\dots x_r}$ has the
critical degree and it makes sense to consider the toric residue
$$\Res_z\left(\frac{z_0\dots z_n}{x_1\dots x_r}\right).$$

As before the monomials $z_0,\dots,z_n$ define a coloring of
$\Sig$ which is now not necessarily disjoint. \rt{main1} implies
the following more general statement.

\begin{Th}\label{T:main2} Let $X$ be an $n$-dimensional projective toric
variety defined by a projective fan $\Sig$. Let $z_0,\dots,z_n$ be
 monomials in the homogeneous coordinate ring
$S=\C[x_1,\dots,x_r]$ such that
\begin{enumerate}
\item $x_1\dots x_r\,|\,z_0\dots z_n$,
\item $z_0,\dots,z_n$ do not vanish simultaneously on $X$.
\end{enumerate}
Then the toric residue
\begin{equation}\nonumber
\Res_z\left(\frac{z_0\dots z_n}{x_1\dots x_r}\right)
\end{equation}
equals the combinatorial degree of the simplicial
coloring of $\Sig$ defined by $z_0,\dots,z_n$.
\end{Th}

\begin{pf} To reduce our theorem to \rt{main1} we choose monomials
$z_0',\dots,z_n'$ such that $z_i'|z_i$ and $z_0'\dots
z_n'=x_1\dots x_r$. Then
\begin{equation}\label{e:transform}
\Res_{z'}(1)=\Res_z\left(\frac{z_0\dots z_n}{x_1\dots
x_r}\right).
\end{equation}
Indeed, the open sets $U_i'=\{z_i'=0\}$, $0\leq
i\leq n$ form a covering $\cal U'$ of $X$. Also the covering $\cU$
by the sets $U_i=\{z_i=0\}$, $0\leq i\leq n$ is a refinement of
$\cU'$. Therefore the cocycles
$$[\omega_{z'}(1)]\in\cZ^n(\cU',\widehat\Omega_X^n)\quad\text{and}\quad
\Big[\omega_{z}\Big(\frac{z_0\dots z_n}
{x_1\dots x_r}\Big)\Big]\in\cZ^n(\cU,\widehat\Omega_X^n)$$
define the same element in $H^n(X,\widehat\Omega_X^n)$.

Now it remains to show that the combinatorial degrees of the two
colorings defined by $z_0,\dots,z_n$ and $z_0',\dots,z_n'$ are the
same. Let $P$ be a polytope whose normal fan is $\Sig$, and
$C=(C_0,\dots,C_n)$ and $C'=(C_0',\dots,C_n')$ be the two
colorings of $P$ defined by $z_0,\dots,z_n$ and $z_0',\dots,z_n'$,
respectively. Consider $G\in\cF(\p P)$ and assume that $G$ belongs
to $C_{i_1}'\cap\dots\cap C_{i_k}'$ with $k$ maximal, where $0\leq
i_1<\dots<i_k\leq n$. Then $G$ belongs to $C_{i_1}\cap\dots\cap
C_{i_k}$ and hence $\psi_C(G)\subset\psi_{C'}(G)$. But this implies that
$\psi_C$ refines $\psi_{C'}$ and so they have the same combinatorial
degree by \rp{subdivision}.

\end{pf}

\begin{Rem} The statement in \re{transform} also follows from the Global
Transformation Law (Theorem~0.1 \cite{CCD}).
\end{Rem}

\section{Monomial ideal}

In this section we
study radical monomial ideals of the homogeneous
coordinate ring $S$ of $X$.
We answer the following question: When does a
homogeneous polynomial $F\in S_\alpha$ of semiample degree $\alpha$ belong to
a radical monomial ideal $I$?

\subsection{Semiample class}\label{S:semiample}
Let $D=\sum_{j=1}^ra_jD_j$ be a $\T$-invariant Cartier divisor on~$X$. Then
it determines a unique continuous function $\phi_D$ linear on each
cone of $\Sig$ such that $\phi_D(v_j)=-a_j$, $1\leq j\leq r$,
where $v_j$ is the first lattice point along the ray
$\rho_j\in\Sig(1)$. Note that if $D$ and $D'$ are linearly
equivalent then $\phi_D-\phi_D'$ is a linear function on $N_\R$.

The divisor $D$ also determines a rational convex polytope $P_D$
in $M_{\R}$:
$$P_D=\{u\in M_{\R}\ |\ \langle u,v_j\rangle\geq-a_j,\ 1\leq j\leq r\}=
\{u\in M_{\R}\ |\ u\geq\phi_D\}.$$
Linearly equivalent divisors have the same polytope up to a parallel
translation.

\begin{Def} A Cartier divisor $D$ is called {\it semiample} if the
corresponding line bundle  $\cO_X(D)$ is generated by global sections.
Also $\alpha=[D]\in A_{n-1}(X)$ is called {\it semiample class}
if $D$ is semiample.
\end{Def}

A divisor $D$ is semiample if and only if the function $\phi_D$ is
convex ([F], Section~3.4).

In \cite{M} A.~Mavlyutov shows that if $\alpha$ is a semiample class
on $X$ then there is a unique toric variety $X'$ and a surjective
morphism $\pi:X\to X'$ such that $\alpha$ is the pull-back of an
ample class $\alpha'$ on $X'$. Here we recall this construction.

Let $\alpha=[D]$ be a semiample class. First,
with the polytope $P_D$ we associate a rational fan $\Sig_D$ as follows.
If $\dim P_D=n$ then $\Sig_D\subset N_\R$ is the normal fan of $P_D$.
In general let $L\subset N_\R$ be the orthogonal complement
to the affine space spanned by~$P_D$. Then for every face $G$ of $P_D$
the cone $\sig_G$ contains $L$ (see \re{normalfan}).
Define~$\Sig_D$ to be the fan in $N_\R/L$ whose cones are the
images $[\sig_G]$ under the quotient map $N_\R\to N_\R/L$.
The fan $\Sig_D$ is rational with respect to the lattice $N/N\cap L$.
Since~$\Sig_D$ is the same for all
representatives $D$ of $\alpha$ we denote it by $\Sig_\alpha$.

The quotient map $N\to N/N\cap L$ extends to a map of fans
$$\tilde\pi:\Sig\to\Sig_\alpha.$$
Indeed, for any maximal cone $\sig\in\Sig$ the restriction of
$\phi_D$ to $\sig$ defines a lattice point $u$ which is a vertex of $P_D$
(by convexity of $\phi_D$). Then $\tilde\pi(\sig)\subset[\sig_u]$. The map
$\tilde\pi$ induces the surjective morphism
$$\pi:X\to X',\quad\text{where}\ \ X'=X(\Sig_\alpha).$$

Note that when $\dim P_D=n$ the map $\tilde\pi:\Sig\to\Sig_\alpha$ is a
refinement. Also $\Sig=\Sig_\alpha$ and $\tilde\pi=id$ when $\alpha$ is
an ample class.

\subsection{Monomial ideal}
Let us now consider a radical monomial ideal $I$ of $S$. We denote by
$z_1,\dots,z_m$ the minimal generating set of $I$.
Similar to the construction of Section~3 these monomials determine
a coloring of rays of
the fan $\Sig$ into $m$ colors: a ray $\rho_j$ has color $i$ if and
only if $x_j$ divides $z_i$. Note that we allow multi-colored rays
and not all the rays of $\Sig$ are colored, in general.

Recall that the {\it irrelevant ideal} of $X$ is the ideal $B$ in $S$
generated by all the products $\prod_{\rho_j\not\in\sig}x_j$ where
$\sig$ runs over all (maximal) cones of $\Sig$ (see \cite{Coxhom}).
Similarly, we define the ideal $B_\alpha$ in $S$ generated by all the
products $\prod_{\tilde\pi(\rho_j)\not\in\sig}x_j$ where $\sig$ runs over all
(maximal) cones of $\Sig_\alpha$ and $\tilde\pi:\Sig\to\Sig_\alpha$
as above. Clearly, $B\subset B_\alpha$ and
$B=B_\alpha$ if $\alpha$ is ample.

\begin{Rem} As mentioned in \cite{CCD}, Section~2
for any ample degree $\alpha$ we have $S_\alpha\subset B$.
This is not true in general, but as
we will see in the proof of the next theorem, $S_\alpha\subset B_\alpha$
when $\alpha$ is semiample.
\end{Rem}

The following is the main result of this section.

\begin{Th}\label{T:monomial} Consider a radical monomial ideal
$I=\langle z_1,\dots,z_m\rangle$ of the homogeneous coordinate ring $S$ of $X$.
Let $\alpha\in A_{n-1}(X)$ be semiample.
The following are equivalent
\begin{enumerate}
\item $S_\alpha\subset I,$
\item $F\in I$, where
$F$ is a generic section of $\alpha$,
\item $B_\alpha\subset I$,
\item no (maximal) cone of $\Sig_\alpha$ contains images of rays of $\Sig$
of all $m$ colors,
\item $\pi(Z_1)\cap\dots\cap\pi(Z_m)=\emptyset$,
where $\pi:X\to X'$ is the morphism determined
by $\tilde\pi:\Sig\to\Sig_\alpha$, and $Z_i$ is the zero
locus of $z_i$ on $X$.
\end{enumerate}

\end{Th}

\begin{pf} (1)$\Leftrightarrow$(2) is clear.

(2)$\Leftrightarrow$(3) Let $D=\sum_{j=1}^ra_jD_j$ be a
representative of $\alpha$ and consider its polytope~$P_D$.
Since $F$ is generic, it is a linear combination
with non-zero complex coefficients of the monomials $\chi^u$ as
$u$ varies over the set $P_D\cap M$ (\cite{F}, Section~3.4).
But $I$ is a monomial ideal thus
$F\in I$ if and only if every $\chi^u$ belongs to
$I$.

In homogeneous coordinates $\chi^u=\prod_{j=1}^rx_j^{\langle
u,v_j\rangle+a_j}$. From this you can see that every
such monomial is divisible by some $\chi^w$, where $w$ is a vertex of $P_D$.
Thus we can assume that $u$ varies over the vertices of $P_D$.
 
When $u$ is a vertex of $P_D$ we have $\langle u,v_j\rangle=-a_j$ if
$\tilde\pi(\rho_j)\in\sig_u$ and $\langle u,v_j\rangle>-a_j$ if
$\tilde\pi(\rho_j)\not\in\sig_u$, where $\sig_u$ is the cone of $\Sig_\alpha$
corresponding to $u$. Therefore we can write
\begin{equation}\label{e:homcoord}
\chi^u=\prod_{\tilde\pi(\rho_j)\not\in\sig_u}x_j^{\langle u,v_j\rangle+a_j},
\end{equation}
where all exponents are positive.
Therefore $\chi^u\in I$ for all vertices
$u$ if and only if
for every maximal cone $\sig$ of $\Sig_\alpha$ there exists $k$ such that
$z_k|\prod_{\tilde\pi(\rho_j)\not\in\sig}x_j$ (by \re{homcoord} and
the fact that in $z_k$ every variable has exponent 0 or 1). By the
definition of $B_\alpha$ the latter is equivalent to
$B_\alpha\subset I$.

(4)$\Leftrightarrow$(5) Suppose there exists a maximal cone
$\sig'$ in $\Sig_\alpha$ which contains the images
$\tilde\pi(\rho_{k_1}),\dots,\tilde\pi(\rho_{k_m})$,
where $\rho_{k_i}$ has color $i$.
It can be readily seen that the morphism $\pi:X\to X'$
maps every orbit $O_{\tau}$, $\tau\in\Sig$ onto the orbit
$O_{\tau'}$, where $\tau'$ is the smallest cone of $\Sig_\alpha$
containing $\tilde\pi(\tau)$, and hence it maps the closure $V(\tau)$ onto
the closure $V(\tau')$. It follows then that for every $0\leq
i\leq n$ the image $\pi(V(\rho_{k_i}))$ contains the closed
orbit $O_{\sig'}$. Since $V(\rho_{k_i})$ is an irreducible
component of $Z_i$ we get a contradiction. The other implication
is similar.

(3)$\Leftrightarrow$(4) Before we give a proof let us take a
closer look at the ideal $B_\alpha$.

We say that a collection of rays $\{\rho_{j_1},\dots,\rho_{j_s}\}$
of $\Sig$ is {\it primitive with respect to $\Sig_\alpha$} if its
image under $\tilde\pi$ does not lie in any of the cones of
$\Sig_\alpha$, but the image of every its proper
subset does. The corresponding collection of
variables $\{x_{j_1},\dots,x_{j_s}\}$ will be also called primitive
with respect to $\Sig_\alpha$.
The following is a slight generalization of the
Batyrev description of the irrelevant ideal of the homogeneous
coordinate ring.

\begin{Lemma} The irreducible components
of the variety $V(B_\alpha)$ of the ideal $B_\alpha$ are the
coordinate planes $V(x_{j_1},\dots,x_{j_s})$ for each collections
$\{x_{j_1},\dots,x_{j_s}\}$ primitive with respect to $\Sig_\alpha$.
\end{Lemma}

\begin{pf} By definition the variety $V(B_\alpha)$ is the
union of $V(x_{\hat\sig}:\sig\in\Sig_\alpha)$, where for every maximal
$\sig\in\Sig_\alpha$ we pick a variable $x_{\hat\sig}$ whose corresponding
ray is mapped outside of~$\sig$.
Since the image of every primitive collection
$\{\rho_{j_1},\dots,\rho_{j_s}\}$ does
not lie in any of the cones of $\Sig_\alpha$ the corresponding collection
$\{x_{j_1},\dots,x_{j_s}\}$ appears among the collections
$\{x_{\hat\sig}:\sig\in\Sig_\alpha\}$. Therefore,
$$\bigcup_{\genfrac{}{}{0pt}{1}{\{\rho_{j_1},\dots,\rho_{j_s}\}}{{ }_{\textrm{primitive}}}}
V(x_{j_1},\dots,x_{j_s})\subseteq V(B_\alpha).$$
On the other hand, every collection $\{x_{\hat\sig}:\sig\in\Sig_\alpha\}$
contains a primitive one. (If a collection of rays
whose image does not lie in any cone is not primitive it contains
a proper subset whose image does not lie in any cone etc.) Thus,
$$V(x_{\hat\sig}:\sig\in\Sig_\alpha)\subseteq V(x_{j_1},\dots,x_{j_s}),$$
for some primitive collection $\{\rho_{j_1},\dots,\rho_{j_s}\}$. The lemma
follows.
\end{pf}

Now we will finish the proof of the theorem. We have
two irreducible decompositions:
$$V(I)=\bigcup_{x_{k_i}|z_i}V(x_{k_1},\dots,x_{k_m})
\quad\text{and}\quad V(B_\alpha)=
\bigcup_{\genfrac{}{}{0pt}{1}{\{\rho_{j_1},\dots,\rho_{j_s}\}}{{ }_{\textrm{primitive}}}}
V(x_{j_1},\dots,x_{j_s})$$
Then $V(I)\subset V(B_\alpha)$ if and only if
for every collection
$\{x_{k_1},\dots,x_{k_m}\}$ with $x_{k_i}|z_i$, there exists a
collection $\{x_{j_1},\dots,x_{j_s}\}$ primitive with
respect to $\Sig_\alpha$ such that
$$V(x_{k_1},\dots,x_{k_m})\subset V(x_{j_1},\dots,x_{j_s}),$$
i.e., $\{\rho_{j_1},\dots,\rho_{j_s}\}\subset\{\rho_{k_1},\dots,\rho_{k_m}\}$.
But the latter is equivalent to (4). Indeed, if $\{\rho_{k_1},\dots,\rho_{k_m}\}$
does not contain any primitive collection then its image lies in some cone
$\sig$ of $\Sig_\alpha$, i.e. $\sig$ contains images of rays of all $m$ colors.
The converse is also clear.
\end{pf}

\section{Homogeneous polynomials of residue one}

Let $X$ be a complete $n$-dimensional toric variety and $F_0,\dots,F_n$ homogeneous polynomials
of degrees $\alpha_0,\dots,\alpha_n$, respectively. Consider the corresponding toric residue
map
$$\Res_F:S_\rho/\langle F_0,\dots,F_n\rangle\to\C.$$
As recently proved by Cox and Dickenstein \cite{CD} this map is an isomorphism if the degrees $\alpha_i$
are semiample and have $n$-dimensional polytopes. Then the computation of the toric residue
reduces to the problem of finding a homogeneous polynomial of residue one. This is an open problem.
The results of the previous sections allow us to say something in this direction.

Assume $X$ is projective defined by a projective fan $\Sig$.
Let $\Sig_i=\Sig_{\alpha_i}$ be the fan associated with
the degree $\alpha_i$ and $\tilde\pi_i:\Sig\to\Sig_i$ the corresponding map as in \rs{semiample}.
We say that a disjoint simplicial coloring $\Sig(1)=\Lambda_0\sqcup\dots\sqcup\Lambda_n$ 
is {\it compatible} with $\Sig_0,\dots,\Sig_n$
if it satisfies the condition of \rt{monomial} for each $\Sig_i$, i.e. no cone of $\Sig_i$
contains images under $\tilde\pi_i$ of rays of all $n+1$ colors, for every $i$.

\begin{Cor} Let $\alpha_0,\dots,\alpha_n$ be semiample degrees on $X$ and $\Sig_0,\dots,\Sig_n$ the corresponding
fans. Fix any disjoint simplicial coloring of $\Sig$ compatible with $\Sig_0,\dots,\Sig_n$:
$$\Sig(1)=\Lambda_0\sqcup\dots\sqcup\Lambda_n.$$
Let $z_i=\prod_{\rho_j\in\Lambda_i}x_j$ be squarefree monomials, $0\leq i\leq n$.
Then for any homogeneous polynomials $F_0,\dots,F_n$ of degrees $\alpha_0,\dots,\alpha_n$
there are homogeneous polynomials $A_{ij}$ such that
$$F_i=A_{i0}z_0+\dots+A_{in}z_n,\quad 0\leq i\leq n.$$
Furthermore,
\begin{equation}\label{e:cor}
\Res_F(\det(A))=\Res_z(1)=\cdeg(\Lambda),
\end{equation}
where $\cdeg(\Lambda)$ is the combinatorial degree of the coloring above. 
\end{Cor}

\begin{pf} The existence of $A_{ij}$ is a consequence of \rt{monomial}. The first statement
in \re{cor} follows from the Global Transformation Law (Theorem 0.1,~\cite{CCD}) and the
second statement follows from \rt{main1}.
\end{pf}

Notice that if the degrees are ample then $\Sig=\Sig_0=\dots=\Sig_n$ and we obtain many different ways of constructing
a polynomial of residue one. In fact, every choice of a disjoint simplicial coloring of $\Sig$ of combinatorial
degree one gives rise to such a polynomial. This extends previously known constructions by D'Andrea and Khetan \cite{AK}
and Cattani, Cox and Dickenstein \cite{CCD}. Also, when the degrees are the same ($\Sig_0=\dots=\Sig_n$) there
are many ways of choosing a compatible coloring of combinatorial degree one. 
It is not true, however, that compatible colorings of combinatorial degree one exist
for any collection of $n+1$ semiample degrees. The simplest counterexample is presented in \rf{counterexample}. 

\begin{figure}[h]
\centerline{
\scalebox{0.6}{
\input{pic1.pstex_t}
             }
           }
         \caption{}
         \label{F:counterexample}
       \end{figure}

Here $P_0$, $P_1$ and $P_2$ are the polytopes of semiample degrees $\alpha_0$, $\alpha_1$ and $\alpha_2$
on the toric surface determined by the complete
fan $\Sig$ on the right. The rays labeled with $i$, for $i=0,1,2$, form the fan $\Sig_i$ corresponding to $\alpha_i$.
One can check that there are no compatible colorings of $\Sig$ of combinatorial degree one. In the forthcoming
paper \cite{KS} further methods are developed to deal with examples like that (in fact, all two-dimensional
examples and a new class of $n$-dimensional examples are considered).


\begin{thebibliography}{GKh2}

\bibitem[{\bf AK}]{AK} C. D'Andrea, A. Khetan,
{\em Macaulay style formulas for toric residues}, preprint
math.AG/0307154.

\bibitem[{\bf BM}]{BM} V. Batyrev, E. Materov,
{\em Toric Residues and Mirror Symmetry}, Moscow Math. J. {\bf 2} (2002), 
no. 3, 435--475.

\bibitem[{\bf CCD}]{CCD} E. Cattani, D. Cox, A. Dickenstein,
{\em Residues in Toric Varieties}, Compositio Math. {\bf 108}
(1997), no. 1, 35--76.

\bibitem[{\bf CaD}]{CaD} E. Cattani, A. Dickenstein,
{\em A global view of residues in the torus}, J. Pure Appl.
Algebra {\bf 117/118} (1997), 119--144.

\bibitem[{\bf CD}]{CD} D. Cox, A. Dickenstein, {\em A Codimension Theorem for Complete Toric Varieties},
preprint, math.AG/0310108.

\bibitem[{\bf CDS}]{CDS} E. Cattani, A. Dickenstein, B. Sturmfels,
{\em Residues and Resultants}, J. Math. Sci. Univ. Tokyo, {\bf 5} (1998), 119--148.

\bibitem[{\bf C1}]{Coxhom} D. A. Cox,
{\em The homogeneous coordinate ring of a toric variety}, J. Algebr. Geom.
{\bf 4} (1995), 17--50.

\bibitem[{\bf C2}]{Coxres} D. A. Cox,
{\em Toric residues}, Arkiv f\" ur Matematik {\bf 34} (1996) 73--96.

\bibitem[{\bf F}]{F} W. Fulton, {\em Introduction to Toric Varieties},
Princeton Univ. Press, Princeton, 1993

\bibitem[{\bf GKh1}]{GKh} O. A. Gelfond and A. G. Khovanskii,
{\em Newton polyhedra and Grothendieck residues}, (in Russian)
Dokl. Akad. Nauk, {\bf 350}, no. 3 (1996), 298--300.

\bibitem[{\bf GKh2}]{GKh2} O. A. Gelfond, A. G. Khovanskii, {\em Toric geometry and
 Grothendieck residues}, Moscow Math.~J. {\bf 2} (2002), no. 1, 99--112.

\bibitem[{\bf KS}]{KS} A. Khetan, I. Soprounov, {\em Partition matrices for polytopes towards
computing toric residue}, preprint math.AG/0406279

\bibitem[{\bf Kh}]{Kh} A. G. Khovanskii,
{\em Newton polyhedra, a new formula for mixed volume, product of
roots of a system of equations}, Fields Inst. Comm., Vol. {\bf 24}
(1999), 325--364.

\bibitem[{\bf M}]{M} A. Mavlyutov, {\em The Hodge structure of semiample
hypersurfaces and a generalization of the monomial-divisor mirror map},
Advances in algebraic geometry motivated by physics, Contemp. Math. 276
(2000).

\bibitem[{\bf S}]{S} I. Soprounov,
{\em Residues and tame symbols on toroidal varieties}, to appear in Compositio Math.,
math.AG/0203114.

\end{thebibliography}
\end{document}